  \documentclass[11pt, oneside]{article}                               

  \usepackage{latexsym}
  \usepackage{amssymb}

  \def\AFOUR{%
  \setlength{\textheight}{9.0in}%
  \setlength{\textwidth}{5.75in}

  \setlength{\topmargin}{-0.375in}%
  \hoffset=-.5in%
  \renewcommand{\baselinestretch}{1.17}%
  \setlength{\parskip}{6pt plus 2pt}%
  }

  \AFOUR                                          

  \makeatletter
  \def\section{\@startsection {section}{1}{\z@}{-3.5ex plus -1ex minus
   -.2ex}{2.3ex plus .2ex}{\large\bf}}
  \def\subsection{\@startsection{subsection}{2}{\z@}{-3.25ex plus -1ex minus
   -.2ex}{1.5ex plus .2ex}{\normalsize\bf}}
  \makeatother

  \makeatletter
  \@addtoreset{equation}{section}
  
  \makeatother

  \newcommand{\nc}{\newcommand}
  \newcommand{\rnc}{\renewcommand}

  \nc{\bea}{\begin{eqnarray}}
  \nc{\eea}{\end{eqnarray}}
  \nc{\be}{\bea}
  \nc{\ee}{\eea}

  \rnc{\a}{\alpha}
  \nc{\ab}{\bar{\a}}
  \nc{\ap}{\a^{+}}
  \nc{\abm}{\ab^{-}}
  \rnc{\b}{\beta}
  \nc{\bb}{\bar{\b}}
  \nc{\bbp}{\bb_{\zb}^{+}}
  \nc{\bm}{\b_{z}^{-}}
  \nc{\oa}{\overline{\a}}
  \nc{\ob}{\overline{\b}}
  \rnc{\gg}{\gamma}
  \rnc{\d}{\delta}
  \nc{\f}{\phi}
  \nc{\fb}{\bar{\phi}}
  \nc{\vf}{\varphi}
  \nc{\p}{\psi}
  
  \rnc{\c}{\chi}
  \nc{\la}{\lambda}
  \nc{\m}{\mu}
  \nc{\n}{\nu}
  \rnc{\o}{\omega}
  \nc{\Om}{\Omega}
  \rnc{\t}{\theta}
  \nc{\eps}{\epsilon}
  \rnc{\S}{\Sigma}
  \nc{\F}{\Phi}

  \nc{\trac}[2]{{\textstyle\frac{#1}{#2}}}

  \nc{\ex}[1]{\mbox{e}^{\,\textstyle#1}}

  \nc{\mat}[4]{\left(\begin{array}{cc}#1&#2\\#3&#4\end{array}\right)}

  \nc{\som}[9]{\left(\begin{array}{ccc}#1&#2&#3\\#4&#5&#6\\#7&#8&#9%
  \end{array}\right)}

  \nc{\tr}{\mathop{\mbox{tr}}\nolimits}
  \nc{\ad}{\mathop{\mbox{ad}}\nolimits}
  \nc{\Tr}{\mathop{\mbox{Tr}}\nolimits}
  \nc{\Det}{\mathop{\mbox{Det}}\nolimits}
  \nc{\rk}{\mathop{\mbox{rk}}\nolimits}
  \nc{\ra}{\rightarrow}
  \nc{\Ra}{\Rightarrow}
  \nc{\LRa}{\Leftrightarrow}
  \nc{\ot}{\otimes}
  \rnc{\ss}{\subset}
  \nc{\nul}{\noindent\underline}
  \nc{\non}{\nonumber\\}

  \nc{\subs}[1]{{\vspace*{0.5cm}}%
  {\noindent\underline{#1}}{\addcontentsline{toc}{subsection}{#1}}%
  {\vspace*{0.3cm}}}

  \nc{\zb}{\bar{z}}
  \rnc{\lg}{\frak{g}} 
  \nc{\lt}{\frak{t}}
  \nc{\lk}{\frak{k}}
  \nc{\lh}{\frak{h}}
  \nc{\pik}{\Pi_{\lk}}
  \nc{\pip}{\Pi_{+}}
  \nc{\pim}{\Pi_{-}}
  \nc{\pih}{\Pi_{\lh}}
  \nc{\jz}{J_{z}}
  \nc{\jzh}{\jz^{\lh}}
  \nc{\jzp}{\jz^{+}}
  \nc{\jzm}{\jz^{-}}
  \nc{\del}{\partial}
  \nc{\dz}{\del_{z}}
  \nc{\dzb}{\del_{\bar{z}}}
  \nc{\az}{A_{z}}
  \nc{\azb}{A_{\bar{z}}}
  \nc{\g}{g^{-1}}
  \nc{\dw}{\Delta_{W}}
  \nc{\Ad}{{\mbox{Ad}}}
  \nc{\ks}{Ka\-za\-ma-\-Su\-zu\-ki}
  \nc{\KS}{\ks}
  \nc{\ksm}{\ks\ model}
  \rnc{\AA}{{\mathbb A}}
  \nc{\BB}{{\mathbb B}}
  \nc{\CC}{{\mathbb C}}
  \nc{\PP}{{\mathbb P}}
  \nc{\cpm}{\CC\PP(m)}
  \nc{\cpn}{\CC\PP(n)}
  \nc{\cp}[1]{\CC\PP(#1)}
  \nc{\gmn}{G(m,m+n)}
  \nc{\gmnk}{\gmn_{k}}
  \nc{\cO}{{\cal O}}
  \nc{\bcO}{\bar{\cO}}
  \nc{\bO}{\bar{O}}
  \nc{\oQ}{\overline{Q}}
  \def \C {\mathbb C }
  \def \R {\mathbb R}

  \def \S {\mathbb S }
  \def \U {1\! \! 1}   
  
  \newtheorem{theorem}{Theorem}[section]
  \newtheorem{proposition}[theorem]{Proposition}
  \newtheorem{lemma}[theorem]{Lemma}
  \newtheorem{corollary}[theorem]{Corollary}

  \newtheorem{remark}[theorem]{Remark}
  \nc{\qed}{\indent \hspace{14.2cm} $\square$}

\begin{document}
 \global\parskip=4pt

 \makeatletter
 \begin{titlepage}
\begin{flushright}
IC/2003/151
\end{flushright}
 \begin{center}
 \vskip .5in
 {\LARGE\bf The Universal Connection and Metrics on Moduli Spaces}\\
\vskip 0.4in
{\bf Fortun\'{e} Massamba$^{a,\, b\, 1}$ and George
  Thompson$^{a}$\footnote{email: massamba/thompson@ictp.trieste.it}}
\end{center}

\hskip 1.8cm a) The Abdus Salam ICTP, 
P.O. Box 586, 34100 Trieste
Italy.

\hskip 1.8cm b) I.M.S.P., B.P. 613, Porto-Novo, B\'{e}nin.

\begin{center}

\vskip .4in
\begin{abstract}
We introduce a class of metrics on gauge theoretic moduli
spaces. These metrics are made out of the universal matrix that
appears in the universal connection construction of M. S. Narasimhan
and S. Ramanan. As an example we construct metrics on the $c_{2}=1$
$SU(2)$ moduli space of instantons on $\mathbb{R}^{4}$ for various
universal matrices.
\end{abstract}

 \end{center}
 \end{titlepage}
 \makeatother

 \begin{small}
 \tableofcontents
 \end{small}
 \setcounter{footnote}{0}

 \section{Introduction}
The aim of this paper is to give, more or less, natural metrics on
gauge theoretic moduli spaces.

The question of alternate metrics to
the standard $L^{2}$ metric arose recently in the context of the AdS/CFT
correspondence. The moduli space in question is the $c_{2}=1$, $SU(2)$
moduli space of instantons on ${\mathbb R}^{4}$. The $L^{2}$ metric is
not `the right one' in that context, essentially because it does not
preserve the conformal invariance inherent in the definition of the moduli
space. A rather remarkable alternative is the information metric which
is built out of $\Tr{F_{A}*F_{A}}$ and its derivatives with respect to
the moduli \cite{H}. The information metric is designed to preserve the
conformal invariance of the theory at hand and yields, for the round
metric on $S^{4}$, the standard Einstein metric on $AdS_{5}$. It is
remarkable in that if one perturbs the metric on $S^{4}$ then to
first order in that perturbation the information metric remains
Einstein \cite{BNT}. There are also a host of other small miracles associated
with this metric.

On other instanton moduli spaces the information metric fails to be a
metric as it becomes highly degenerate. The reason for this is that it
is overly gauge invariant, meaning that it is invariant under gauge
transformations that depend on the moduli. For example a convenient
parameterization of the $c_{2}=1$, $SU(n)$ instanton
moduli space on $\mathbb{R}^{4}$ is the one where the $SU(2)$
instanton is embedded in $SU(n)$ and then one acts with rigid $SU(n)$
gauge transformations to obtain the general instanton. But all the
rigid gauge transformations leave $\Tr{F_{A}*F_{A}}$ invariant and so
drop right out of the information metric leaving us with the $SU(2)$
parameters only. Rather more dramatically one sees that the
information metric is in fact zero on moduli spaces of flat or
parabolic bundles. So, while the $L^{2}$ metric does not keep some of
the properties we would like, it does have the advantage of being a metric! 

The problem then, as set out at the start of this introduction is to
find other natural metrics. The
$L^{2}$ metric is made out of the gauge connection directly while the
information metric is constructed from the curvature 2-form. These are
natural objects in the theory and that is why these metrics are also
natural. These, however, do not exhaust the natural objects that are
available to us. 

In 1961 M.S. Narasimhan and Ramanan \cite{NR} introduced the concept
of a universal connection. This connection plays the same role as that
of the universal bundle construction for bundles. Specifically any
connection on the Principal bundle of interest can be obtained by
pull-back from the universal connection (the bundle itself is obtained
by pull-back of the universal bundle). The important conclusion there
is that any $U(n)$ connection\footnote{There is no restriction on the
structure group, we have fixed on $U(n)$ here for ease of
presentation.} can be expressed as 
\be
A = U^{\dagger}d\, U
\ee
where $U$ is an $m\times n$ rectangular matrix satisfying
\be
U^{\dagger}\, . \, U = {\mathbb I}_{n}.
\ee
Since any such connection is made up of $n^{2}d$ real `functions' on a
d-dimensional manifold it is apparent that there is a raw lower bound
on $m$ namely that $2m \geq n(d+1)$. This bound is very difficult to
meet. Indeed the bound one gets depends precisely on how the matrices
$U$ are constructed.

In any case, the observation of M.S. Narasimhan
and Ramanan was turned into a
powerful tool for self-dual connections on 4-manifolds in the ADHM
construction. That construction amounts to a method for obtaining a
from a universal connection the required matrix $U$.

Infact the universal connection appears in the construction of many
moduli spaces. We take the attitude that the $U$ matrices are also
natural objects in the gauge theory and so should be used in the
construction of metrics. 
Our construction gives us metrics on moduli spaces where the
parameterization of the moduli space is contained in the
connection. See our main assumption in the next section. When there
are `matter' or `Higgs' fields present one should use that data as well
in the construction of metrics but this depends on the details of the
equations one is trying to solve so we do not enter into this, except
loosely in the Conclusions. As an example
one can check that for instantons on flat ${\mathbb R}^{4}$ with $U$ given
by the ADHM construction one of the proposed metrics is the
non-degenerate $AdS_{5}$ metric.

The paper is organized as follows. We delay describing
the universal connection construction, and so how one may arrive at
the matrices $U$, till Section \ref{NRsect}. Consequently, we require
the readers indulgence until then and ask that they take on faith results that
will be proven in that section. In the
next section we introduce possible metrics on the space of
connections and ultimately on the moduli space of interest. In Section
\ref{NRsect} the universal connection construction is finally
given. An improvement in the Abelian case is also presented. Metrics
on the $c_{2}=1$ $SU(2)$ instanton moduli space are considered in
detail in Section \ref{ims} by way of example of the general
construction. Finally, in the Conclusions, we end with many open questions.

\section{Universal Metrics}\label{sect1}

The title of this section is perhaps misleading. We mean that these are
metrics made from the matrices $U$. Suppose that we have a parametrized
family of connections, with parameters $t^{i}$ such that
\be
A(t) = U^{\dagger}(t) d U(t).\label{pu}
\ee
Denote the derivative $\partial / \partial t^{i}$ by
$\partial_{i}$. We denote by ${\mathcal G}$ the space of maps from $M$
to $G$. We also denote by ${\mathcal G}_{t}$ the space of maps from $T
\times M$ to $G$, where $T$ is the space of moduli. So infact for
fixed $t^{i} \in T$, $g(t,x) \in {\mathcal G}$.

{\bf Assumptions:} We demand that the parameters are `honest' parameters,
that is, we demand that the parameterization is complete. Furthermore,
we will presume that we are working on a smooth part of the moduli
space. The two assumptions imply that there are no non-zero vectors $v^{i}(t)$
such that $v^{i}(t)\partial_{i} A(t) = 0$ at any point $t^{i}$ of the
moduli space under consideration. We will refer to these assumptions
as the main assumption.

Introduce the projector,
\be
P = U.U^{\dagger}
\ee
which is an $m \times m$ hermitian matrix, satisfying
$$
P^{2} = P.
$$
The projector, $P$, is not only gauge invariant since $U\rightarrow
U.g^{\dagger}$ for $g \in {\mathcal G}$ but also invariant under
${\mathcal G}_{t}$.
We have the following simple,
\begin{lemma}\label{pdl}{{\rm
If $\partial_{i} P = 0$ 
then $\partial_{i} A = d_{A} A_{i}$; where 
$A_{i} = U^{\dag}\, \partial_{i} U$ and the covariant derivative
is $d_{A}= d +[ A, $ .
}}
\end{lemma}
{\bf Proof:}
This is by direct computation. $\partial_{i} P = 0$ implies that
$$
\partial_{i} U^{\dag} = - U^{\dag}\,\partial_{i}U \, U^{\dag}, \;\;\; 
\partial_{i} U = - U\,\partial_{i}U^{\dag}\,U.
$$
For the variation of
the connection we have, thanks to these two equations,
\bea
\partial_{i}A & = & \partial_{i}U^{\dagger}. dU + U^{\dagger} d
\partial_{i}  U \nonumber \\
& = &  -\,U^{\dag}\partial_{i} U \, A + A \,U^{\dag}\,\partial_{i} U 
+ d\left(U^{\dag}\, \partial_{i} U\right)\nonumber\\
&=&d_{A} A_{i}\nonumber 
\eea
\hskip 14.3cm $\square$

One can consider the connection form to be a form in a higher
dimension, that is on $T\times M$, with the components in the $T$
direction being $A_{i}= U^{\dagger} \partial_{i} U$. With this
understood Lemma \ref{pdl} says that $v^{i} \partial_{i}P=0$ only if
$v^{i}F_{i\, \mu}=0$ ($F_{i\, \mu}$ are the mixed components of the
curvature 2-form on $T\times M$).

We now introduce some universal metrics. 
Set
\bea
g^{0}_{ij} & = & \int_{M} \Tr{(\partial_{i} P * \partial_{j} P)}
\label{0metric}. \\
g^{1}_{ij} & = &  -\int_{M} \Tr A_{i} \, * \, A_{j} ,\label{1metric}
\eea
where $*$ is the Hodge star operator. $g^{0}$ is invariant under
${\mathcal G}_{t}$ while $g^{1}$ is only invariant under
${\mathcal G}$. The metrics that one gets naturally descend to
${\mathcal A}/{\mathcal G}$ and finally to the moduli space.

\begin{theorem}\label{Main}{\rm{Let $M$ be a compact closed manifold and $U$ a
      universal matrix for some family of connections then there is a linear
      combination of the components of the quadratic forms $g^{0}$
      (\ref{0metric}) and $g^{1}$(\ref{1metric}) which is a metric on the
      moduli space.}}
\end{theorem}
{\bf Proof:} There are two cases to consider:

First suppose that there
is no vector $v^{i}$ such that 
$v^{i}\partial_{i}P=0$, then $g^{0}$ is a metric. This follows from
the fact that $*$ and $\Tr$ (on Hermitian matrices) are positive
definite so that for $g^{0}$
to be degenerate, that is for $g_{ij}v^{i}\, v^{j}=0$ for some $v^{i}$,
we must have $v^{i}\partial_{i}P=0$.

Secondly suppose that there is a vector $v^{i}$ such that
$v^{i}\partial_{i}P=0$ then $g^{0}$ is degenerate in this
direction. However, by Lemma \ref{pdl} we have that
$v^{i}\partial_{i}A = d_{A} \, v^{i} A_{i}$ and $v^{i}A_{i}$ cannot be
zero as that would contradict our main assumption. Positive definiteness of
$*$ and negative definiteness of $\Tr$ (on anti-Hermitian matrices)
guarantee that 
$g^{1}_{ij}v^{i}\, v^{j} \neq 0$. 

Consequently we can always organize
for some linear combination of the components of $g^{0}$ and $g^{1}$ to yield
a non-degenerate symmetric quadratic form on $T$.

\qed

We have a kind of converse to the theorem,
\begin{corollary}[to Theorem \ref{Main}]\label{cor1}{\rm{Given the conditions of
      the theorem if $A_{i}=0$ then $g^{0}$ is a metric on the moduli space.}}
\end{corollary}
{\bf Proof:} If $A_{i}=0$ then there is no $v^{i}$ such that
$v^{i}\partial_{i}P=0$ since if there was we would conclude
$v^{i}\partial_{i}A=0$ contradicting our main assumption.

\qed

To see that we can apply Corollary \ref{cor1} of Theorem \ref{Main}
directly we quote a 
\begin{lemma}\label{main}{\rm{If $U$ is the NR (M.S. Narasimhan,
      Ramanan) matrix then
$$
A_{i}=U^{\dagger} \partial_{i} \, U = 0.
$$
}}
\end{lemma}
{\bf Proof:} This is delayed till Section \ref{NRsect} where all the
definitions will also be available.

As the construction in \cite{NR} applies to any connection we learn,
from Corollary \ref{cor1}, that any moduli space can be given the
metric $g^{0}$ providing the universal matrix used is a NR matrix.

If one wishes to make use of $U$ matrices other than those which are
NR matrices then (\ref{0metric}) can fail to be, but need not fail to be, a
metric. In the proof of Theorem \ref{Main} we saw that degeneracy of
$g^{0}$ only comes from having connections whose dependence on some moduli
is through gauge transformations. This is precisely the situation that
we described in the Introduction for the data for instantons for
higher rank and which plagues the information metric.

In cases of this type one has
that
some of the moduli, say $s^{a}$, are obtained by a gauge transformation
on a connection $A_{0}$ which depends on moduli $r^{\alpha}$, that is
$A(s,r) = (U_{0}(r)h(s,r))^{\dagger}d(U_{0}(r)h(s,r))=
h^{\dagger}(s,r) A_{0}(r) 
h(s,r ) +  h^{\dagger}(s,r)dh(s,r) $ then $\partial_{a}A(s,r) = d_{A}
A_{a}(s,r)$ where $A_{0}(r)= U^{\dagger}_{0}(r)dU_{0}(r)$ and 
$A_{a}(s,r)= h^{\dagger}(s,r)\partial_{a} h(s,r)$. Furthermore, we have that
$P$ depends on all the coordinates $r^{\alpha}$ but $\partial_{a}P =0$.
Now $g_{ab}^{0}=g_{a\alpha}^{0}=0$ however,
$$
g^{1}_{ab} = -\int_{M} \Tr {\left( h^{\dagger}(s,r)\partial_{a}h(s,r) *
  h^{\dagger}(s,r)\partial_{b}h(s,r)\right)}.
$$
Thus, in this situation, one can consider as a metric a linear
combination of $g^{0}$ and $g^{1}$.
\begin{proposition}\label{g1g2}{\rm{Let $M$ be a compact, closed
      manifold and $U$ a universal matrix for a family of connections
      with $U=U_{0}(r)\, h(s,r)$ where $g^{0}$ for $U_{0}$ is a metric
      on the part of the moduli space parameterized by $r^{\alpha}$ and $h\in
      \mathcal{G}_{s}$ as above, then a linear combination of $g^{0}$
      and $g^{1}_{ab}$ is a metric on the moduli space.}}
\end{proposition}
\qed

\begin{remark}{\rm{Notice that we are not saying that $A_{\alpha}=
g^{\dagger}(s,r)\partial g/ \partial r^{\alpha}$ is zero. Such a
condition is not required for $g^{0}$ to be a metric.}}
\end{remark}

\subsection{More Metrics from the Universal Connection}

One of the problems we are faced with is that for non-compact
manifolds the integrals that go into defining $g^{0}$ and $g^{1}$ may
well not converge. To improve the situation we add a `damping'
factor. This generalization gives us many possible metrics even in the
compact case. 

Let
\be
\Phi(U) = *^{-1} \Tr( dP * dP)
\ee 
So far we
had only used a volume form on $M$, however $\Phi(U)$ requires a
metric. Note that $\Phi(U)$ is invariant under
$\mathcal{G}_{t}$. 

As an aside note that the mass dimension of $\Phi(U)$ is 2 which means
that it is like a mass term for the gauge field
$A_{\mu}$. Infact it is gauge invariant albeit highly
non-local and non-polynomial (in the gauge field). One has
$$
\int_{\mathbb{R}^{4}} d^{4}x \; \Phi(U) = \int_{\mathbb{R}^{4}}
d^{4}x\, 
\Tr\, \left( - A_{\mu}\, A^{\mu} + 
\partial_{\mu}U^{\dagger} \, \partial^{\mu}U \right).
$$
which is much more suggestive of a mass for the gauge field.

Set
\bea
& & g^{0,\, \alpha}_{ij} = \int_{M} \Phi(U)^{\alpha}\, \Tr ( \partial_{i}P*
\partial_{j}P), \label{g0a}\\
& & g^{1,\, \beta}_{ij} =
-\int_{M}\Phi(U)^{\beta} 
\Tr {\left( U^{\dagger}\partial_{i}U * 
  U^{\dagger}\partial_{j}U\right)}\label{g1b}.
\eea

\begin{proposition}\label{3prop}{\rm{$g^{0,\, \alpha}$ and $g^{1,\,
	\beta}$ clearly have the following properties.
\begin{enumerate}
\item They are gauge invariant (that is under $\mathcal{G}$).
\item For $\alpha = \beta = d/2$ where $\dim M= d$, the metric on $M$
  enters only through its conformal class.
\item For $M$ non-compact with $\alpha$ and $\beta$ suitably large the
  integrals in (\ref{g0a}, \ref{g1b}) formally converge provided that
  $\Phi(U)$ has some suitable integrability properties.

\end{enumerate}
}}\end{proposition}
\qed

As far as the third point of the proposition is concerned it may be
natural to suppose that $\Phi(U)$ be integrable (as indicated by the
analogy with a mass term). However, weaker integrability conditions
may also suffice in certain cases.

\section{The Universal Connection}\label{NRsect}

M. S. Narasimhan and Ramanan \cite{NR}
prove that, at least locally, any connection, $A$ on a bundle
$\mathcal{P}$ can be expressed as
\be
A = i U^{\dagger} d U \label{univ}.
\ee

The set up is as follows \cite{NR}. One begins with the Stiefel manifold,
$V(m,n)$ of all unitary $n$-frames through the origin of ${\mathbb
  C}^{m}$, thought of as a $U(n)$ principal bundle over the Grassman
manifold $G(m,n) \equiv U(m)/\left(U(m) \times U(m-n)\right)$. There
is a connection on this bundle which has the following
description. Denote the components of the $n$-planes by $v_{i} =
\sum_{a=1}^{m}  e_{a}\, S_{ai}$ where the $(e_{a})$ are a canonical
basis for ${\mathbb C}^{m}$ and $i = 1, \dots ,\, n$. By orthogonality
one requires that $S^{\dagger}\, . \, S = {\U}_{n}$ (the unit
$n\times n$ matrix). Denote the
matrix valued function which 
associates to each $n$-plane its matrix $S_{ai}$ by the same
letter. Then 
\be
S^{\dagger}\, d\, S = \omega
\ee
is a canonical connection on the bundle. This is the universal connection. The
main theorem of M.S. Narasimhan and S. Ramanan is that any connection
on a principle $U(n)$-bundle, $\mathcal{P}$ over a $d$-dimensional
manifold $X$, is obtained by
pullback of $\omega$ from a differentiable bundle homomorphism from
$\mathcal{P}$ to $V(m,n)$ for some sufficiently large $m$.

The theorem tells us that any connection may be expressed in the form
(\ref{univ}) for an $m\times n$ matrix $U$ providing that $m$ is
sufficiently large. They provide a
lower bound on $m$, 
$m \geq (d+1)(2d+1)n^{3}$ will do, but it is
not a very efficient one. We will see below that for the local problem for
$U(1)$ bundles one can infact do much better than requiring $m=(d+1)(2d+1)$.

\subsection{The Narasimhan-Ramanan Matrix}

 In their paper \cite{NR} Narasimhan and Ramanan not only prove the
 existence of the universal connection but also give a construction of
 the matrices $U$. The way they do this is to pass from a local
 construction of $U$ to a more global one. They certainly do not give
 the most minimal form of $U$ but, nevertheless, their procedure is
 the only one we know of that will produce the required matrix $U$ for
 any connection.

 \begin{lemma}{{\rm \cite{NR} 
 Let $V$ be an open subset of $\R^{d}$ and $W$ a relatively compact open
 subset whose closure is contained in $V$. For every differential form
 $\alpha$ of degree 1 on $V$ with values in $\mathfrak{u}(n)$ (the space of
 skew-Hermitian matrices), there exist differentiable functions
 $\phi_1,\cdots\phi_{m'}$ in $W$ with values in the space
 $\mathcal{M}_{n}(\C)$ of $\left(n\times n\right)$ complex matrices such that 
 \begin{enumerate}
 \item $\sum_{j=1}^{m'}\phi_j^{\ast}\phi_j=\U_{n}$,  and
 \item $A= \sum_{j=1}^{m'}\phi_j^{\ast}d\phi_j$.
 \end{enumerate}
where $m'=(2d+1)n^{2}$.
 }}
 \end{lemma}

The Narasimhan-Ramanan matrices are of a very particular form. Write 
\be
A = i\sum_{\mu=1}^{d}\sum_{r=1}^{n^{2}}\lambda_{r,\mu}f_r dx_{\mu}
\ee
where $\left(f_1\cdots\cdots f_{n^2}\right)$ a set of positive define
matrices which form a base for the complex Hermitian matrices over
the reals, such that $||f_r||=1$ for every $r$, here $||f||$ being the
norm as a  linear transformation. According to the proof the lemma,
there exists $p_{r,\mu}.\;\;q_{r,\mu}$ and $h_r$ strictly positive
differentiable functions such that
\be
\lambda_{r,\mu}&=&p_{r,\mu}^2-q_{r,\mu}^2
, \nonumber\\
h_r^2(x)&=&\frac{1}{n^2}\U_{n}-\left[\sum_{\mu=1}^{d}(p_{r,\mu}^2+
  q_{r,\mu}^2) \right]f_r
\ee

The matrix $U$ that Narasimhan and Ramanan propose is,
\be
 U(A)=\left[\begin{array}{c}
 \Phi_1\\
 \Phi_2\\
 \Phi_3
 \end{array}\right],
\ee
where $\Phi_1$ has $d(n \times n)$ components defined by the functions
 $p_{r,\mu}e^{ix_{\mu}} \cdot g_{r}$, $\Phi_2$ has $d(n
 \times n)$ components defined by the functions $q_{r,\mu}\, e^{-ix_{\mu}}\cdot
 g_{r}$, $\Phi_3$ has $(n \times n)$ functions defined by
 $h_r$ and $g_{r}$ is the positive square root of $f_{r}$. Think of
 all of these as vectors with entries $(n \times n)$ matrices.

There are two, somewhat surprising, results that can be deduced from
the lemma and its proof.
\begin{proposition}\label{prop1}{\rm{In general the Narasimhan-Ramanan
      matrix for the gauge
      transform of a connection is not the gauge transformation of the
      Narasimhan-Ramanan matrix
      of the original connection, i.e. it is not equivariant,
\be
U(A^{g}) \neq U(A) . g \label{cov}
\ee}}
\end{proposition}
{\bf Lemma \ref{main}} Let $A(t)$ be a family of connections
      parameterized by $t^{i}$ and let $U(t)$ be the corresponding family
      of Narasimhan-Ramanan universal matrices. Then,
\be
U^{\dagger}(A) \partial_{i} U(A) = 0. \label{inv}
\ee
where $\partial_{i}=\partial/\partial t^{i}$.

The first proposition is evident from the construction of the matrix
$U(A)$ (one should refer to \cite{NR} for the details of that
construction).

{\bf Proof of Lemma \ref{main}:}

 From the definitions we have,
 $$
 \Phi_{1}^{\dagger}\,
 \partial_{i}\Phi_{1}  = \frac{1}{2}\sum_{r, \mu} \partial_{i}p_{r,
   \mu}^{2}\, f_{r}, \;\;\; \Phi_{2}^{\dagger}\,
 \partial_{i}\Phi_{2}  = \frac{1}{2}\sum_{r, \mu} \partial_{i}q_{r,
   \mu}^{2}\, f_{r},\;\;\; \Phi_{3}^{\dagger}\,
 \partial_{i}\Phi_{3}  = \frac{1}{2}\sum_{r} \partial_{i} h_{r}^{2},
 $$
 so that
 $$
 \sum_{a=1}^{3}\Phi_{a}^{\dagger}\,
 \partial_{i}\Phi_{a} = \partial_{i}\frac{1}{2}\left(
 \sum_{r,\mu}(p_{r, \mu}^{2}+q _{r, \mu}^{2})f_{r} + \sum_{r} h_{r}^{2}
 \right)  =  \partial_{i}\frac{1}{4n^{2}}\U_{n \times n} =
 0.\;\;\;\;\; 
 $$
\qed

From Proposition \ref{prop1} we learn that there are universal
parameterizations which are not gauge covariant. In particular, it is
difficult to see how to construct invariants from the
Narasimhan-Ramanan matrix apart from those made out of the curvature
2-form. This suggests that, in this case, one should already work on a
slice of the space of connections with this parameterization.

Somewhat more mysterious is Lemma \ref{main}. It implies a
special case of Proposition \ref{prop1}. Suppose that $U$ is some
parameterization matrix, not necessarily the NR
matrix, which satisfies (\ref{inv}) and is covariant $U(A^{g})=
U(A).g$ for all $g \in \mathcal{G}_{t}$
then we get into a contradiction. Inserting the covariance condition into
(\ref{inv}) we find that
$$
g^{\dagger}\partial_{i}\, g = 0,
$$
which implies that $g$ cannot be an arbitrary gauge
transformation but rather only one that does not depend on the
parameters $t^{i}$.

\subsection{Abelian Universal Connections}

The NR matrices are very messy to deal with, as can be seen
from the details of their construction. However, in the Abelian case
one can simplify the discussion and construction somewhat.

\begin{lemma}{
{\rm
There exists fixed real-valued functions 
$r_{i}$, for 
$i=1, \dots , d-1$, on $\R^{d}$ such that 
any $U(1)$ connection $A$, which is
pure gauge at infinity, can be expressed as 
\be
A=-\sum_{i=1}^{d-1}\theta_{i}\, d r^2_{i} + d\bar{\theta}.
\ee
for some functions $\overline{\theta}$ and $\theta_{i}$, with $i=1,
\dots, d-1$. }}
\end{lemma} 

{\bf Proof:} Let $\partial_{j} r^2_{i}=\delta_{ij}\partial_{i}
r^2_{i}$ for $i=1, \dots, d-1$ and such this derivative does not
vanish anywhere except at infinity, i.e. $\partial_{i} r^2_{i}
\rightarrow 0$ as $ |x| \rightarrow \infty $. This means that we can
invert $\partial_{i} r^2_{i}$ everywhere. Let $A$ match $d
\overline{\theta}$ at infinity. Set,
$$
\theta_{i}=-\frac{1}{\partial_{i} r^2_{i}}
\left(A_{i}-\partial_{i}\bar{\theta}_{d}\right)
$$
and in any case we have that $A_{d}= \partial_{d}
\overline{\theta}$, from which we can solve for
$\overline{\theta}$. Plugging back in establishes the lemma. 

A set of functions that have the required property are $
r^2_{i}=c_{i}(\exp(x_i)+1)^{-1}$ where $c_{i}$ are constants.

\qed

\begin{corollary}\label{cor3}{{\rm A Universal matrix for any $U(1)$ connection
      $A$ on $\R^{d}$, which is pure gauge at infinity, is given by 
\begin{eqnarray*}
U=\left[\begin{array}{l}
 r_1 e^{-i\left(\theta_1+\theta_{d}\right)}\\
 r_2 e^{-i\left(\theta_2+\theta_{d}\right)}\\
\vdots\\
r_{d-1} e^{-i\left(\theta_{d-1}+\theta_{d}\right)}\\
r_{d} e^{-i\theta_{d}}
\end{array}\right],
\end{eqnarray*}
where $r_{d}^{2}=1 -
\sum_{i=1}^{d-1}r_{i}^{2}$.
and $\theta_{d}=\bar{\theta}-\sum_{i=1}^{d-1} r^2_{i}\theta_{i}$.}}
\end{corollary}

{\bf Proof:} We have that 
\be
A=iU^{\dag} dU=\sum_{i=1}^{d-1} r^2_{i} d(\theta_{i}+\theta_{d}) +
d\theta_{d}=-\sum_{i=1}^{d-1}  dr^2_{i} \, \theta_{i} +
d\overline{\theta} \nonumber
\ee
and by the lemma this is the general form for any such connection.

\qed

The NR matrices, in the $U(1)$ case, are $2d+1 \times 1$ matrices. The
connection matrix above is a $d \times 1$ matrix and since a
connection form is essentially given by $d$ functions we may consider
this to be an optimal parameterization.

\section{Instanton Moduli Space}\label{ims}

As already mentioned a convenient parameterization of the $c_{2}=1$,
$SU(n)$ instanton moduli space on $\mathbb{R}^{4}$ is the one where the $SU(2)$
instanton is embedded in $SU(n)$ and then one acts with rigid $SU(n)$
gauge transformations to obtain the general instanton. From previous
sections we know that the metric coming from these gauge
transformations will be caught by $g^{1}$ and the $SU(2)$ part will
come from $g^{0}$. In this section we concentrate on the $SU(2)$ part.

The moduli space $\mathcal{M}_{SU(2)}^{1}$ is known to be
 five-dimensional, with the topology of the open 5-ball, thus
 parametrized by five
 moduli, namely four coordinates $a^{\mu}$ (the centre of the
 instanton) and one scale $\rho$. The parameter $\rho$ measures the
 size of the instanton, and zero size corresponds to delta function
 or so called singular instantons. $\rho=0$ corresponds to the boundary
 of the 5-ball, $S^{4}$, where the moduli $a^{\mu}$ are the coordinates
 on the $S^{4}$.

The following sections are relatively brief as all the details are
essentially computational and many steps are skipped.

\subsection{The ADHM Universal Matrix and its Metric}

The ADHM construction gives the universal connection form for the
 instanton that we are interested in \cite{ADHM}. These authors find a
 $4 \times 2$ representation for the rectangular matrices required to
 parameterize an instanton, namely the explicit expression of $U$ for
 the self-dual gauge potential is
 \begin{eqnarray}
  U(x)&=&\frac{1}{\sqrt{(x-a)^2+\rho^2}}\left[\begin{array}{c}
  (x-a)^{\mu}\bar{\sigma}_{\mu}\\\\
  -\rho\U_{2\times2}
  \end{array}\right]\label{U}
  \end{eqnarray}
 We use the notation
 \be
 \sigma_{\mu}= (\U_{2\times 2}, i \tau_{a}), \;\;\;
  \overline{\sigma}_{\mu}= (\U_{2\times 2},-  i \tau_{a}),
 \ee
 with $\tau_{a}$ the usual Pauli matrices.

A straightforward calculation leads to
\bea
 \Tr{\left(\partial_{a^{\mu}
 } P\, 
 \partial_{a^{\nu}}P\right)} & = & \frac{4
   \delta_{\mu
 \nu}\,  \rho^{2}}{\left( (x-a)^{2} + \rho^{2}\right)^{2}}\nonumber \\
 \Tr{\left(\partial_{\rho
 }P\, \partial_{a^{\nu}}P \right)} & =&\frac{4
   \rho\,  (x-a)_{
 \nu}}{\left( (x-a)^{2} + \rho^{2}\right)^{2}}\nonumber \\
  \Tr{\left(\partial_{\rho
 }P\,  \partial_{\rho}P \right)}  & = &\frac{4
    (x-a)^{2}}{\left( (x-a)^{2} + \rho^{2}\right)^{2}} .
 \eea
We do not have to work to calculate $\Phi(U)$. By translational
invariance one can replace derivatives with respect to $x^{\mu}$ with
derivatives with respect to $-a^{\mu}$ so that we find
 \be
 \Tr\,\Phi(U)= \frac{16  \rho^{2}}{\left( (x-a)^{2} + \rho^{2}\right)^{2}}
 \ee

Our next task is to determine the metric $g^{0, \, \alpha}$. For
$\alpha >1/2 $ the integrals 
converge and we consider $\alpha$ in this range. By
rotational invariance the integral defining $g^{0,\, \alpha}_{\rho
\, a^{\mu}}$ is zero. Likewise, by translational invariance, the other
integrals do not depend on $a^{\mu}$. The metric on the moduli space
is, therefore, of the form
\be
ds^{2} = \rho^{1-\alpha}\, A(\alpha)\left( d\vec{a}^{2} + B(\alpha)d\rho^{2}
\right) \label{gmetric}
 \ee
where the dependence on $\rho$ is determined by dimensional
arguments. The coefficients $A(\alpha)$ and $B(\alpha)$ are determined
on doing the integrals and they are both non-zero, thus this is a
metric. When $\alpha=1$ this is proportional to the $L^{2}$
metric. When $\alpha = 4/2$, $ds^{2}$ is proportional to the
$\mathrm{AdS}_{5}$ metric and this had to be so as a consequence of
Proposition \ref{3prop}. Thus with $\alpha = 2$ the universal metric
exhibits the nice feature of the information metric namely that it is
Einstein. 

For completeness we list the values of the coefficients
\bea
 A(\alpha) &=& (4)^{2\alpha +1}\,\pi^{2}\, \left[2\alpha-1\right]
\frac{\Gamma\left(2\alpha-1\right)}{\Gamma\left(2(\alpha+1)\right)}
\nonumber  \\
 B(\alpha) &=&\frac{2}{2\alpha-1}.\nonumber
 \eea

\subsection{The NR Universal Matrix and its Metric}

The construction of the NR matrix is somewhat involved, so here we
will outline some of the ingredients and leave the details for
\cite{Massamba}. 

For the instanton moduli space
$\mathcal{M}^1_{SU(2)}$, the explicit expression for the connection of
interest, $A_{\mu}$, is given by 
\begin{eqnarray}\label{c}
 A_{\mu}=\eta_{\mu\nu}^a\frac{(x-a)_{\nu}}{\left(x-a\right)^2+\rho^2}\tau_a,
 \end{eqnarray}
where $\eta_{\mu\nu}^a$ are the 't Hooft eta-symbols, a basis for
self-dual two-form on $\R^4$, and $\tau_{a}$ are the Pauli
matrices. The basis of positive definite Hermitian matrices $f_{r}$ are chosen
to be $f_{i}= (\tau_{i} + 3 \mathbb{I}_{2} )/4$ for $i=1,\; 2,\; 3$
and $f_{4}= \mathbb{I}_{2}$. The $\lambda_{r,\mu}$ are given by the expressions
\be
 \lambda_{i,\mu}&=&4\eta_{\mu\nu}^{i}\frac{(x-a)_{\nu}}{ \left(x-a
   \right)^2 +\rho^2  }, \;\; i=1,\; 2, \; 3\\
 \lambda_{4,\mu}&=&- \frac{3}{4} \sum_{i=1}^{3} \lambda_{i,\mu}.
\eea

One can see directly from the construction of $p_{r\mu}$ and
$q_{r\mu}$, as spelled out in the proof of the Lemma on page 565 of
\cite{NR}, that these functions depend on the position of the
instanton, $a^{\mu}$, only through the combination $(x-a)^{\mu}$. This
is due to the fact that the $a_{r, \mu}$ (that appear on page 566),
for the single $SU(2)$ instanton, are functions of the scale $\rho$
and not of $a^{\mu}$. Consequently, once more by dimensional arguments,
the $g^{0\, \alpha}$ metric will take the form given in
(\ref{gmetric}) and only the coefficients $A(\alpha)$ and $B(\alpha)$
need to be determined.

This agreement of $g^{0}$, up to change of parameters, of the ADHM and
NR universal matrices is not a `generic' situation. The proof of the Lemma just
cited holds if the $a_{r, \mu}$ are replaced by $a_{r, \mu} + k_{r,
  \mu}$ where the $k_{r,\mu}$ are positive functions. For example one
could consider $k_{r,\mu}= c_{r , \mu} \exp{(- d_{r,\mu} a^{2}/
  \rho^{2})}$ with $c_{r , \mu}$ and $d_{r,\mu}$ positive
constants satisfying suitable conditions so that $h_{r}$ is well
defined. In this case the universal matrix would lead to a metric 
with a highly non-trivial dependence on $a^{2}$.

\section{Conclusions}

There are many more metrics that one can form from the universal
matrices $U$. For example, both the
$L^{2}$ and information metrics can be written in terms $U$. One can
also construct other `damping' factors. If we set
$$
\phi = *^{-1} dP*dP
$$
then powers of terms of the form 
$$
\Tr \phi^{\alpha_{1}}\dots \Tr \phi^{\alpha_{k}}
$$
with $\alpha_{i}$ positive integers, give us more positive definite
damping factors. Proposition \ref{3prop} holds if one replaces $\Phi$
with the damping factors we have just introduced.

One question that has not been addressed is: how are metrics that come
from different parameterizations, meaning different universal
matrices, related? In some cases are they related simply by a
diffeomorphism? We do not have an answer to this series of
questions. However, there 
is a set of cases in this direction where we have a simple
\begin{proposition}{\rm{Let $U$ be a
 $m\times n$ universal matrix and $M$ a $k\times n$ universal matrix
 with $k\ge m$ so that 
 $$
 A_{\mu}=iU^{\dagger}\partial_{\mu}U =iM^{\dagger}\partial_{\mu}M.
 $$
If
$$
  M =\frac{1}{\sqrt{N}}\left[\begin{array}{c}
  U\\
  U\\
  \vdots\\
  U
 \end{array}\right] ,
$$
with $k= Nm$, then the universal metrics of $U$ and $M$ agree.}}
\end{proposition}
\qed

We have also not investigated the geometry of the metrics that have
been introduced. This appears to require a more detailed knowledge of
the moduli space that one is dealing with as, for example, of the
instanton moduli space of the previous section.

Many gauge theory moduli spaces also involve `matter' fields and
consequently any metric on the moduli space would probably need to
involve those objects as well. Suppose $\Psi$ is such a field, that is
a section of some associated bundle to the Principal bundle (tensored
with other bundles). Think of $\Psi$ as being valued in some tensor
product representation $V
\otimes \dots \otimes V$ where $V$ is the $n$-dimensional
representation of $SU(n)$. In this case the covariant derivative
$$
d_{A} \Psi = U^{\dagger} \otimes \dots \otimes U^{\dagger}\, d \, ( U
\otimes \dots \otimes U \, \Psi)
$$
and one can form gauge invariant
combinations $U
\otimes \dots \otimes U \,  \Psi$ and from this construct gauge
invariant terms to be added to the universal metric.

Finally, we would also like to know if there exists a construction
which produces equivariant universal matrices. We note that the $U$
matrices in the Abelian case in Corollary \ref{cor3} are indeed equivariant.

{\bf Acknowledgments:} It is a pleasure to thank M. Blau, K. Narain,
M.S. Narasimhan and T. Ramadas for many useful
discussions. F. Massamba would like to thank the Abdus Salam ICTP for
a fellowship. This research was supported in part by EEC contract
HPRN-CT-2000-00148.

\rnc{\Large}{\normalsize}

\end{document}